\documentclass[final,3p]{elsarticle}
\usepackage{amsmath}
\usepackage{amssymb}
\usepackage{graphicx}
\usepackage{snapshot}
\usepackage{units}
\usepackage{esint}
\usepackage{subfigure}
\usepackage{lineno}

\journal{Applied Mathematics Letters}

\begin{document}

\begin{frontmatter}

\title{A Polynomial Approximation for Arbitrary Functions}

\author[mike]{Michael A. Cohen}
\ead{mike@cns.bu.edu}
\author[can1,can2]{Can Ozan Tan\corref{cor1}}
\ead{cotan@partners.org}

\cortext[cor1]{Corresponding author. SW052, Cardiovascular Research Laboratory, Spaulding Rehabilitation Hospital Cambridge, 1575 Cambridge Street, Cambridge, MA 02138. Tel: +1-617-758-5510. Fax: +1-617-758-5514}

\address[mike]{Department of Cognitive and Neural Systems, Boston University, Boston, MA}
\address[can1]{Department of Physical Medicine and Rehabilitation, Harvard Medical School, Boston, MA}
\address[can2]{Cardiovascular Research Laboratory, Spaulding Rehabilitation Hospital, Cambridge, MA}

\begin{abstract}
We describe an expansion of Legendre polynomials, analogous to the
Taylor expansion, to approximate arbitrary functions. We show that
the polynomial coefficients in Legendre expansion, thus, the whole
series, converge to zero much more rapidly compared to the Taylor
expansion of the same order. Furthermore, using numerical
analysis with sixth-order polynomial expansion, 
we demonstrate that the Legendre polynomial approximation yields
an error at least an order of magnitude smaller than the analogous
Taylor series approximation. This strongly suggests that Legendre
expansions, instead of Taylor expansions, should be used when global
accuracy is important.
\end{abstract}


\end{frontmatter}

\linenumbers

\section{Introduction}

The question whether and how a given function can be expressed approximately
by polynomials are of great importance in theory as well as in practice.
For example, by definition, an explicit, finite formula is unavailable
for transcendental functions, and instead, an appropriate polynomial
approximation is chosen to replace the function. Because
polynomials, particularly the ones of low order, are easy to manipulate,
this approach provides computational speed with minimal penalty in
accuracy.

A natural candidate for polynomial approximation is a truncated Taylor
expansion, typically at the midpoint of the interval where the 
approximation is most accurate. Taylor's theorem and the Weierstrass
approximation theorem \cite{weierstrass1985} asserts the possibility
of local approximation of an arbitrary function $y(x)$. Moreover, 
the approximation accuracy improves as the degree of the polynomial increases.
However, this improvement comes at the expense of complexity and computational
speed. This expense can be substantially reduced if the 
function can be approximated to the same accuracy with a lower degree
polynomial.

Here, we show analytically that an arbitrary function can be approximated via 
Legendre polynomials using \emph{non-uniformly} spaced points on an interval 
as the input, and that at least for some functions, 
approximation with Legendre polynomials yields a substantially higher
accuracy and faster convergence compared to Taylor expansion of the same 
order (i.e., with the same number of non--zero coefficients). We further 
demonstrate the improvement in accuracy over Taylor expansion numerically,
using the sine, exponential, and entropy functions. 

\section{Theory}


Consider the problem of estimating the instantaneous slope of the
curve mapping the the output of the function $y(x)$ to $x\in[a,b]$. The
formula for the slope in a linear regression $y(x)\sim x+\epsilon$
for uniformly spaced continuous points $x$ over the interval $[a,b]$
is given by 

\begin{equation}
s=\dfrac{\dfrac{1}{b-a}\intop_{a}^{b}E\left(y(x)-E(y(x))\right)(x-E(x)))dx}{\dfrac{1}{b-a}\intop_{a}^{b}(x-E(x))^{2}dx}\label{eq:slope}
\end{equation}

\noindent where $E(x)$ denotes the expectation of $x$. 
Because $x$ is uniform in $[a,b]$, the denominator
of the equation~\ref{eq:slope} can be written as

\begin{equation}
\frac{1}{b-a}\int_{a}^{b}(x-E(x))^{2}dx=\sigma^{2}_{[a,b]}=\dfrac{(b-a)^{2}}{12}\label{eq:denom}
\end{equation}

\noindent where $\sigma^2_{[a,b]}$ is the variance in the interval $[a,b]$. The numerator of the equation~\ref{eq:slope} can be written as

\begin{eqnarray}
\dfrac{1}{b-a}\intop_{a}^{b}E\left(y(x)-E(y(x)\right)(x-E(x)))dx & = & \dfrac{1}{b-a}\intop_{a}^{b}E(y(x)-E(y(x)))xdx\nonumber \\
 & = & \dfrac{1}{b-a}\left(\intop_{a}^{b}y(x)xdx-\dfrac{1}{2}\dfrac{b^{2}-a^{2}}{(b-a)}\intop_{a}^{b}y(x)dy(x)\right) \nonumber\\
 & = & \dfrac{1}{b-a}\left(\left[\intop_{a}^{b}y(x)xdx-\dfrac{b+a}{2}\intop_{a}^{b}y(x)dy(x)\right]\right)\label{eq:intp0} \\
 & = & \dfrac{1}{b-a}\left[\intop_{a}^{b}y(x)\left(x-\dfrac{b+a}{2}\right)dx\right]\label{eq:intp1}
\end{eqnarray}

\noindent which can be solved using integration by parts:

\begin{eqnarray}
\intop_{a}^{b}y(x)xdx & = & \dfrac{b^{2}y(b)-a^{2}y(a)}{2}-\intop_{a}^{b}\dfrac{x^{2}}{2}y'(x)dx \label{eq:intpp0}\\
\intop_{a}^{b}y(x)dy(x) & = & by(b)-ay(a)-\intop_{a}^{b}xy'(x)dx\label{eq:intpp1}
\end{eqnarray}

\noindent where $y'=\frac{dy(x)}{dx}$. Plugging equations \ref{eq:intpp0} and \ref{eq:intpp1}
into equation~\ref{eq:intp0}, equation~\ref{eq:intp1} can be written as

\begin{eqnarray*}
\dfrac{1}{b-a}\left[\intop_{a}^{b}y(x)\left(x-\dfrac{b+a}{2}\right)dx\right] & = & \frac{1}{b-a}\left[\dfrac{b^{2}y(b)-a^{2}y(a)}{2}-\intop_{a}^{b}\dfrac{x^{2}}{2}y'(x)dx\right.\nonumber \\
 &  & \left.-\dfrac{b+a}{2}\left(by(b)-ay(a)-\intop_{a}^{b}xy'(x)dx\right)\right]\nonumber \\
 & = & \frac{1}{b-a}\left[-\dfrac{ab}{2}\left(y(b)-y(a)\right)+\dfrac{b+a}{2}\intop_{a}^{b}xy'(x)dx-\intop_{a}^{b}\dfrac{x^{2}}{2}y'(x)dx\right]\nonumber \\
 & = & \frac{1}{b-a}\left[-\dfrac{ab}{2}\intop_{a}^{b}y'(x)dx+\dfrac{b+a}{2}\intop_{a}^{b}xy'(x)dx-\intop_{a}^{b}\dfrac{x^{2}}{2}y'(x)dx\right]\nonumber 
\end{eqnarray*}

\begin{eqnarray}
 & = & \dfrac{1}{b-a}\left[\intop_{a}^{b}y'(x)\left(-\dfrac{ab}{2}+\dfrac{b+a}{2}x-\dfrac{x^{2}}{2}\right)dx\right]\nonumber \\
 & = & \dfrac{1}{b-a}\left[\intop_{a}^{b}y'(x)\left(-\dfrac{ab}{2}-\dfrac{1}{2}\left(x-\dfrac{b+a}{2}\right)^{2}+\dfrac{(a+b)^{2}}{8}\right)dx\right]\nonumber \\
 & = & \dfrac{1}{b-a}\left[\intop_{a}^{b}y'(x)\left(-\dfrac{1}{2}\left(x-\dfrac{b+a}{2}\right)^{2}+\dfrac{(b-a)^{2}}{8}\right)dx\right]\label{eq:intpp2}
\end{eqnarray}

\noindent Finally, combining equations~\ref{eq:intpp2} and \ref{eq:denom},
equation~\ref{eq:slope} can be rewritten as

\begin{equation}
s=\dfrac{12}{(b-a)^{2}}\left[\intop_{a}^{b}y'(x)\left(-\dfrac{1}{2}\left(x-\dfrac{b+a}{2}\right)^{2}+\dfrac{(b-a)^{2}}{8}\right)dx\right]\label{eq:slope-2}
\end{equation}

\noindent which is just an average with respect to a quadratic kernel
that is centered at the midpoint of the interval $[a,b]$ and zero
at the ends. Equation~\ref{eq:slope-2} allows estimation of the
instantaneous slope over not just the points that are uniformly spaced,
but all points in the interval $[a,b]$.

This result for estimation of the slope is far more general. It provides
a least squares polynomial approximation to an arbitrary function
on the interval $[a,b]$. To see this, consider the shifted Legendre
polynomials of order $n$, defined by Rodrigues' Formula \cite{abramowitz1964handbook,gradshteyn1988tables}:

\begin{eqnarray}
P_{n,[b,a]}(x) & = & \dfrac{1}{2^{n}n!}\dfrac{d^{n}}{dx^{n}}\left(1-\left[\dfrac{2x-a-b}{b-a}\right]^2\right)^{n}\nonumber \\
 & = & \dfrac{1}{n!(b-a)^{n}}\dfrac{d^{n}}{dx^{n}}\left[\left(b-x\right)(x-a)\right]^{n} \label{eq:rodrigues}
\end{eqnarray}

\noindent which are orthogonal functions with respect to the $L_2$ inner product 

\begin{equation}
<P_{n,[b,a]},P_{m,[b,a]}>=\dfrac{b-a}{2n+1}\delta_{nm}\label{eq:innerprod}
\end{equation}

\noindent where $\delta_{nm}$ denotes the Kronecker delta, equal to $1$ if $m = n$ 
and to $0$ otherwise. Furthermore, Legendre polynomials of order $1$ to $n$
are the same as the orthogonal polynomials obtained by the Gram-Schmidt
process on the polynomials $\{1, x, x^2, \dots\}$ with respect to the inner product 
given by equation~\ref{eq:innerprod} up to a constant multiplication factor \cite{abramowitz1964handbook,gradshteyn1988tables}. Therefore, by adding the basis functions $1,x,\ldots,x^{n}$, we obtain
the $n^{\mathrm{th}}$--order polynomial fit to an arbitrary function
$y$ on the interval $[a,b]$ as

\begin{equation}
y=\sum_{n=0}^{\infty}\dfrac{<y,P_{n,[b,a]}>}{<P_{n,[b,a]},P_{n,[b,a]}>}P_{n,[b,a]}\label{eq:polyfit}
\end{equation}

\noindent By completing the sum of squares, the expected integrated squared
residual error for equation~\ref{eq:polyfit} can be written as

\begin{equation}
\dfrac{1}{b-a}\left(\left\Vert y\right\Vert ^{2}-\sum_{n=0}^{\infty}\dfrac{<y,P_{n,[b,a]}>^{2}}{\left\|P_{n,[b,a]}\right\|^{2}}\right)\label{eq:isre}
\end{equation}

The term $\dfrac{<y,P_{n,[b,a]}>}{<P_{n,[b,a]},P_{n,[b,a]}>}$ in
equation~\ref{eq:polyfit} has a simple and telling interpretation.
Note that 

\begin{eqnarray}
<y,P_{n,[b,a]}> & = & \dfrac{1}{n!(b-a)^{n}}\intop_{a}^{b}y\dfrac{d^{n}}{dx^{n}}\left[\left(b-x\right)(x-a)\right]^{n}dx\label{eq:legendre_approx} \\
 & = & \dfrac{(-1)^{n}}{n!(b-a)^{n}}\intop_{a}^{b}\dfrac{d^{n}y}{dx^{n}}\left[\left(b-x\right)(x-a)\right]^{n}dx\label{eq:legendre}
\end{eqnarray}

\noindent which follows from equation~\ref{eq:rodrigues}. The integral 
in equation~\ref{eq:legendre} can be solved using integration by
parts, and none of the boundary terms appear in the solution because
$\dfrac{d^{k}}{d^{k}x}(x-a)^{n}(x-b)^{n}=0$ for $x=a,b$ and $k=0,\cdots n-1$.
Moreover,

\begin{equation}
<y,P_{n,[b,a]}> = \dfrac{(b-a)^{n}(-1)^{n}n!}{(2n+1)!}E_{\beta(n+1,n+1;a,b)}\left[\dfrac{d^{n}y}{dx^{n}}\right]
\end{equation}

\noindent where $\beta(a,b;c,d)$ is a beta probability distribution
with $(a,b)$ degrees of freedom, shifted and scaled so that it lies
on the interval $[c,d]$, and $E_{\beta(n+1,n+1;a,b)}(\cdot)$ is the 
expected value with respect to this distribution. Therefore, the 
contribution of each of the $n^\mathrm{th}-$order Legendre polynomial 
to the residual error for equation~\ref{eq:polyfit} (shown in 
equation~\ref{eq:isre}) can be rewritten as

\begin{equation}
\dfrac{<y,P_{n,[b,a]}>^{2}}{\left\|P_{n,[b,a]}\right\|^{2}}=\left[\dfrac{(b-a)^{2n-1}}{(2n+1)!{2n \choose n}}\right]\left(E_{\beta(n+1,n+1;a,b)}\left[\dfrac{d^{n}y}{dx^{n}}\right]\right)^{2}\label{eq:legendre-2}
\end{equation}

Equation~\ref{eq:legendre-2} is analogous
to the Taylor expansion, and uses a scaled value of the expected value
of the $n^{\mathrm{th}}$ derivative in the expansion. Note that the distribution 
$\beta(n+1,n+1;a,b)$ is a parabola raised to the $n^{\mathrm{th}}$ power and scaled to
integrate to 1, so it rapidly approximates $\left(x-\dfrac{a+b}{2}\right)$.
The mean and the variance of the distribution $\Omega\sim\beta(a,b;c,d)$
is given by 
\begin{eqnarray}
E(\Omega) & = & \dfrac{a+b}{2}\\
V(\Omega) & = & \dfrac{1}{2N+3}\dfrac{\left(b-a\right)^{2}}{4}
\end{eqnarray}

\noindent and the entire series (equation~\ref{eq:polyfit}) takes the form

\begin{equation}
y=\sum_{n=0}^{\infty}\left[\dfrac{(-1){}^{n}E_{\beta(n+1,n+1;a,b)}\left[\dfrac{d^{n}y}{dx^{n}}\right]}{n!{2n \choose n}}(b-a)^{n-1}\right]P_{n,[b,a]}\label{eq:legendre_series}
\end{equation}


The series in equation~\ref{eq:legendre_series} converges much more
rapidly than the Taylor series because the coefficient ${(-1)^{n}}/{{2n \choose n}n!}$
converges to zero much more rapidly than the analogous term ${1}/{n!}$
in Taylor series. Furthermore, if $\left|{d^{n}y}/{dx^{n}}\right|\leq K\alpha^{n}$
for arbitrary constants $K$ and $\alpha$, it is easy to show (by Stirling\textquoteright{}s
formula and by the nature of Legendre functions) that the squared norm
of the $n^{\mathrm{th}}$ term is bounded above by

\begin{equation}
\dfrac{K^{2}\alpha^{2n}(b-a)^{2n-1}n!}{[2n]^{2}!}^{2}\approxeq\dfrac{K^{2}\alpha^{2n}e^{2n}(b-a)^{2n-1}}{16^{n}n^{4n}}=\dfrac{K^{2}}{(b-a)}\left(\dfrac{\alpha^{2}e^{2}(b-a)^{2}}{16n^{4}}\right)^{n}\label{eq:conv_bound}
\end{equation}

Note that the multiplier of the numerator, but not that of the denominator, of
equation~\ref{eq:conv_bound} is bounded. Let $T_{n_0}$ denote the first term of 
the series in equation~\ref{eq:legendre_series}. If we choose its order, $n_0$, 
such that $n_{0}^{4}\geq\mbox{\ensuremath{{\alpha^{2}e^{2}(b-a)^{2}}/{8}}}$,
then the squared norm of each successive term $\ensuremath{T_{k+n_{0}}}$ has a norm 

\begin{equation}
2^{-k}\left(\dfrac{a^{2}e^{2}}{16(n_{0}+k)}\right)^{n_{0}+k}<2^{-k}T_{n_{0}} 
\end{equation}

\noindent for $n_{0}+k\geq n_{0}$. This shows that the sum of squares of the 
terms after the $n_{0}^{\mathrm{th}}$
term is at most twice the size of the previous term, and the series
is dominated by a geometric sequence with a decay
constant $\unitfrac{1}{2}$. Because one can start with any arbitrary 
$n\geq n_{0}$, the upper bound of the absolute value of $n^\mathrm{th}$ term 
(equation~\ref{eq:conv_bound}) can be reduced arbitrarily by simply increasing
the number of terms, $n$. 

Demonstrating convergence in absolute value is also simple because
the Legendre polynomials are normalized so that they are bounded by
one. The bound, equivalent to that in equation~\ref{eq:conv_bound}, for
the maximum of the $n^{\mathrm{th}}$ term is given as: 

\begin{equation}
\dfrac{K\alpha^{n}(b-a)^{n-1}n!}{[2n]!}\approxeq\dfrac{K\alpha^{n}(b-a)^{n-1}e^{n}}{\sqrt{2}(4)^{n}n^{2n}}=\dfrac{K}{\sqrt{2}}\left(\dfrac{\alpha(b-a)e}{(4)n^{2}}\right)^{n}
\end{equation}

\noindent This bound also indicates that for the same $n_{0}$ (as
above), the maximum norm of each successive term is again less than $2^{-k}\max(T_{n_0})$. 
Similar to above, for large $n_{0}$, the absolute value of
the terms $n>n_{0}$ is less than twice the absolute value of the
term $n_{0}$, which can be arbitrarily small. Thus, the Legendre Polynomial expansion 
converges at a very rapid rate.

\section{Numerical Results \label{sec:Numerical}}

\begin{figure}[h]
 \subfigure[Entropy function]{
   \includegraphics[width=0.33\textwidth,angle=0]{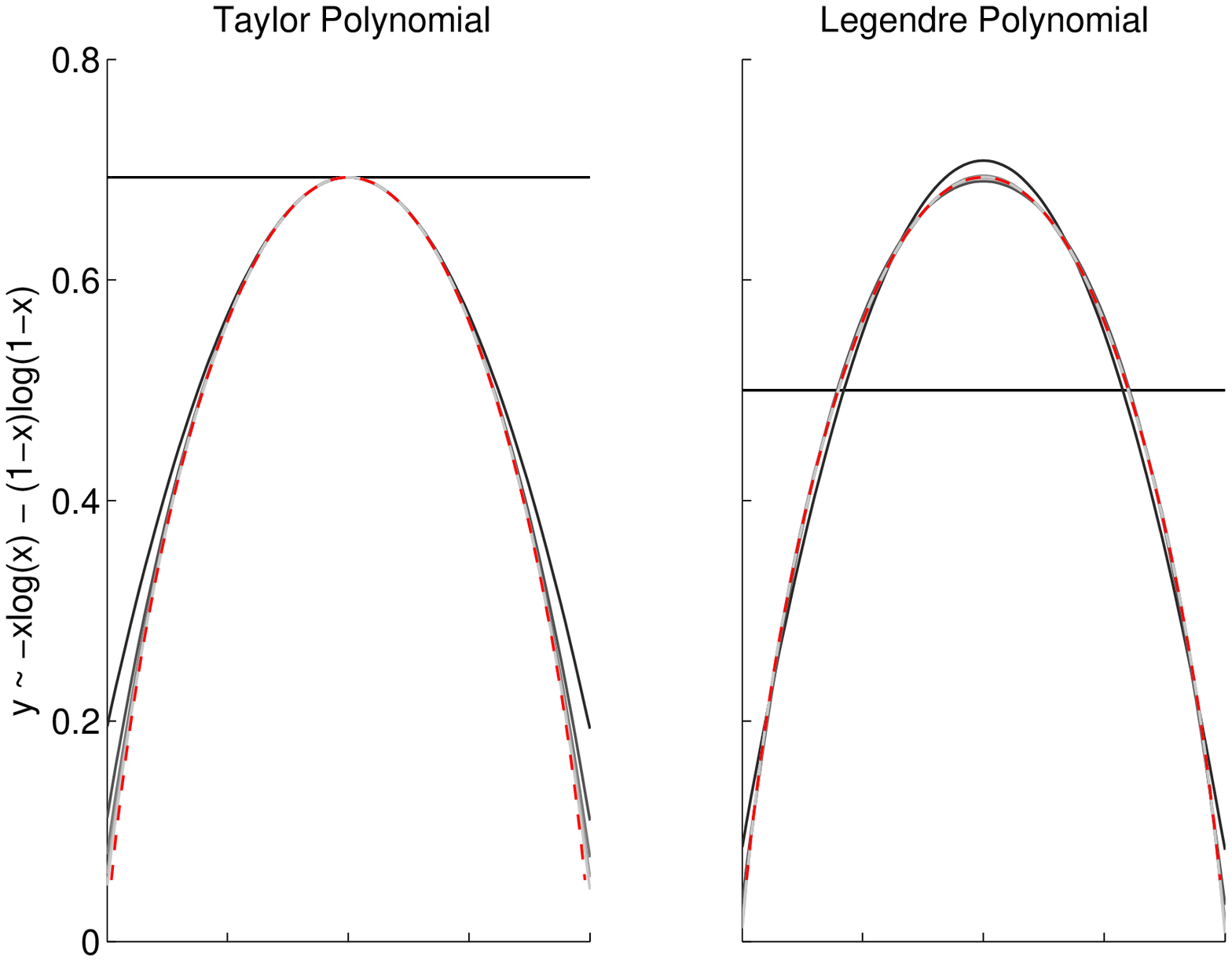}\label{fig:subFit1}
 }
 \subfigure[Sine function]{
   \includegraphics[width=0.33\textwidth,angle=0]{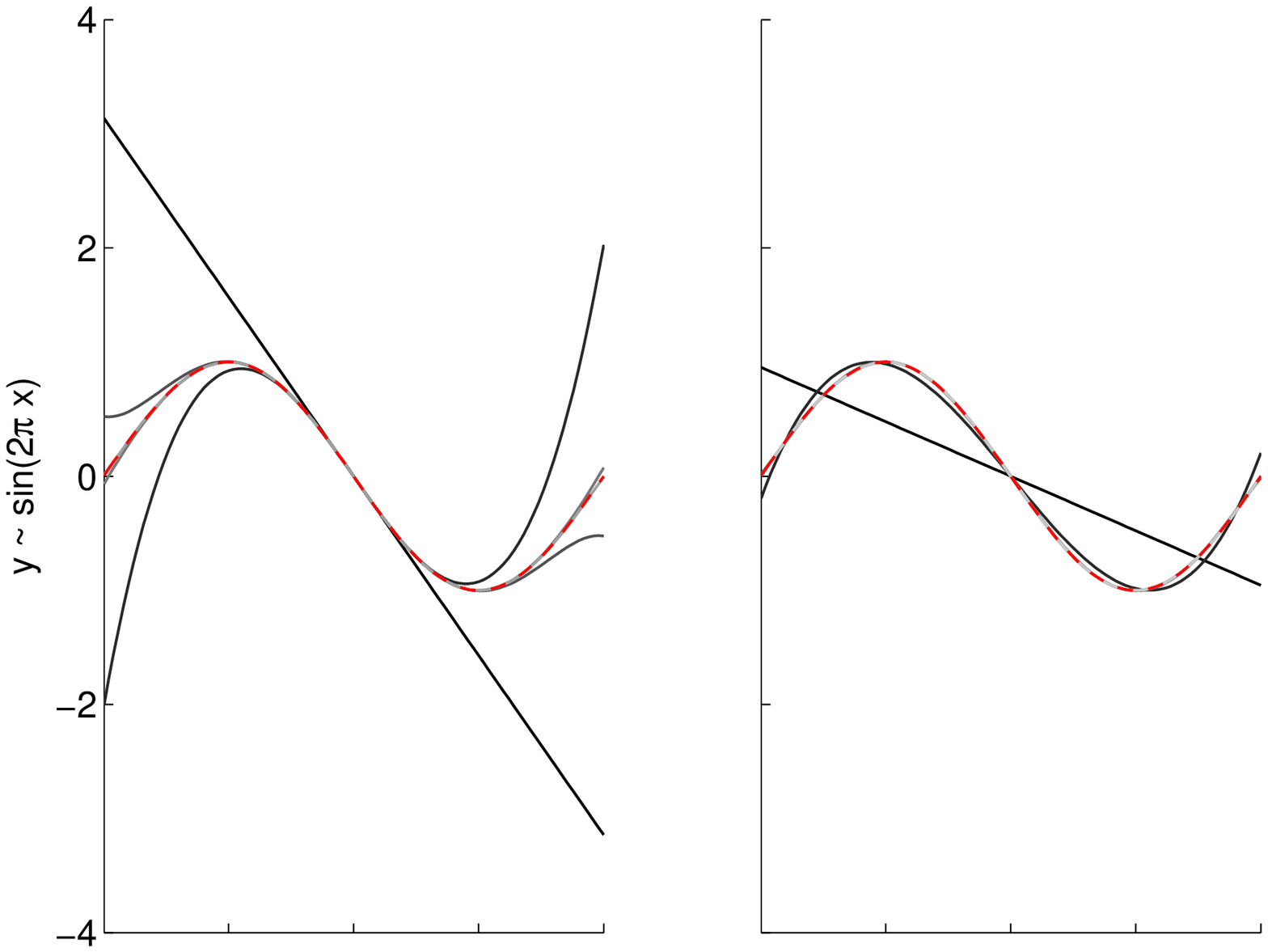}\label{fig:subFit2}
 }
 \subfigure[Exponential function]{
   \includegraphics[width=0.33\textwidth,angle=0]{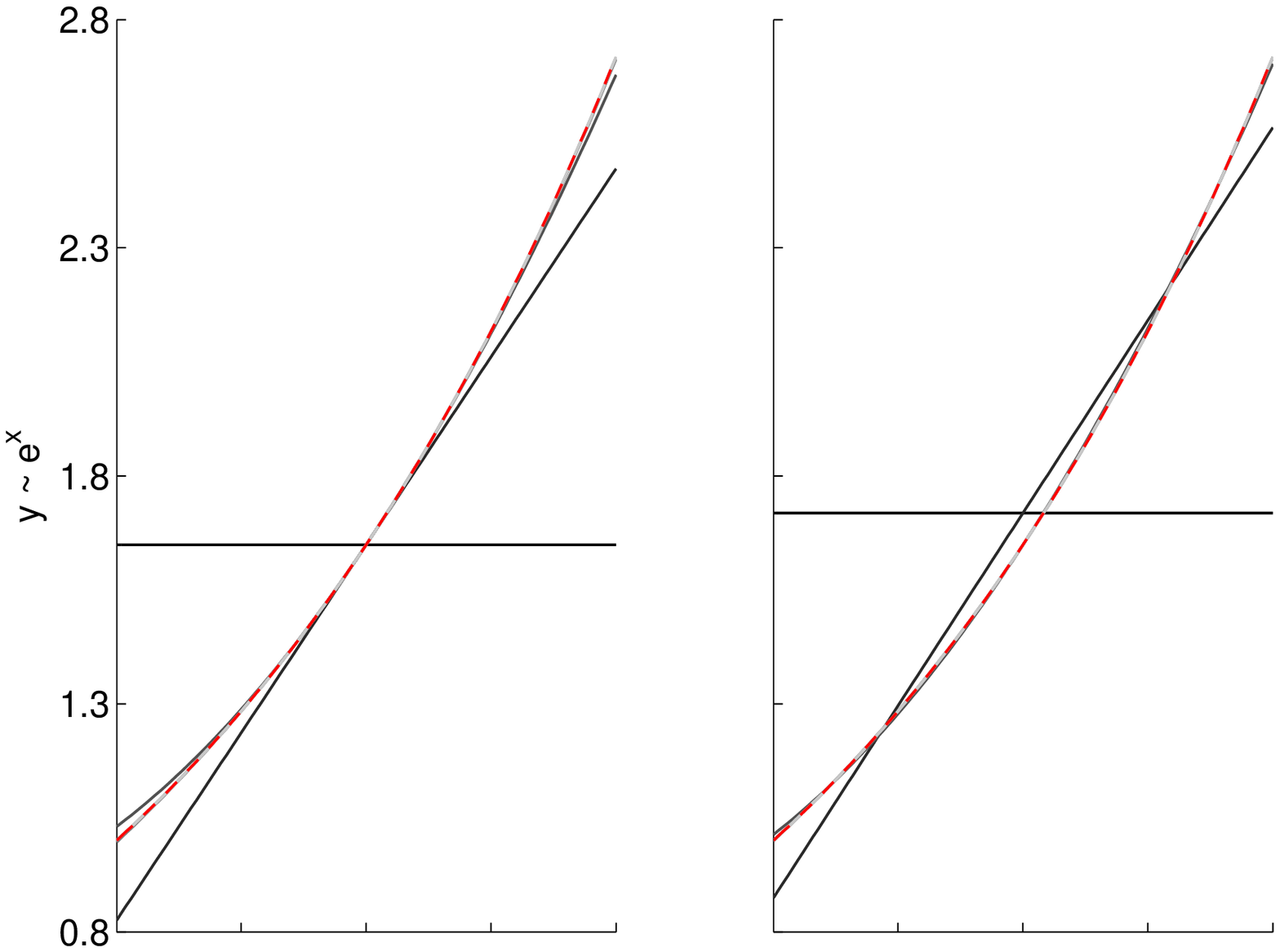}\label{fig:subFit3}
 }
 \caption{\label{fig:TayvsLagFit}Taylor and Legendre polynomial approximations
  at $x \in [0,1]$ Dashed red traces show the actual function evaluated for 
  each $x$. Values of polynomial approximations are shown with solid traces, 
  and successively darker colors indicate approximations using one (light-gray) 
  to six (black) non--zero polynomial coefficients.}
\end{figure}

There is no simple computable remainder (corresponding to the one
in Taylor series expansion) in Legendre polynomial approximation, 
but only bounds (see above).
Therefore, we demonstrate the improvement in accuracy
via an empirical approach using three representative
transcendental functions: the sine function ($\sin(2\pi x)$), the
exponential function ($\exp(x)$), and the entropy of a Bernoulli
trial ($H(p)=-p\log p-(1-p)\log(1-p)$) in the interval $[0,1]$.
Note that all three functions satisfy the boundedness condition 
$\left|{d^{n}y}/{dx^{n}}\right|\leq K\alpha^{n}$
above. For comparison, we approximate
the functions with Taylor and Legendre polynomials using
at most six non-zero coefficients (Figure~\ref{fig:TayvsLagFit}).
Because the sine function has odd symmetry about $x=1/2$ and the
entropy function $H(p)$ is odd, the expansions of these function
involve, respectively, polynomial degrees of 11 and 10 of even and
odd order. We estimate the signal-to-noise ratio (SNR,
in decibels) for approximations as

\begin{equation}
\mathrm{SNR}(f,\hat{f})=10\log_{10}\dfrac{\intop_{0}^{1}f^{2}}{\intop_{0}^{1}(f-\hat{f})^{2}}
\end{equation}

\noindent to quantify the accuracy in approximation (Table~\ref{tab:impr}).

\begin{table}[h]
\caption{Signal-to-noise ratios (in decibels) for the Taylor (T) and Legendre
(L) polynomial approximations. Note the higher SNR for the Legendre
polynomial approximation. The accuracy of Legendre polynomial approximation 
is significantly greater than that of Taylor approximation of equal order.
\label{tab:impr}}
\begin{center}
\begin{tabular}{c|cccccc}
 & \multicolumn{2}{c}{$H(p)$} & \multicolumn{2}{c}{$\sin(2\pi x)$} & \multicolumn{2}{c}{$\exp(x)$}\tabularnewline
Order  & T  & L  & T  & L  & T  & L\tabularnewline
\hline 
1  & 11.12  & 11.20  & -5.53  & 4.0  & 5.95  & 9.11\tabularnewline
2  & 25.54  & 29.09  & .96  & 20.56  & 20.76  & 29.78\tabularnewline
3  & 42.64  & 50.60  & 14.23  & 44.32  & 27.74  & 41.13\tabularnewline
4  & 61.84  & 74.64  & 32.18  & 73.10  & 32.29  & 42.04\tabularnewline
5  & 82.73  & 100.64  & 53.76  & 105.73  & 35.67  & 55.12\tabularnewline
6  & 105.06  & 128.22  & 78.31  & 141.49  & 38.36  & 60.05\tabularnewline
\end{tabular}
\end{center}
\end{table}

Consistent with the theoretical considerations in above,
these results show a rapid convergence of Legendre polynomials 
for increasing number of polynomial coefficients, and an improved
accuracy compared to Taylor polynomials of the same order. 
Note that using six coefficients, even the least improvement
(in the case of the exponential function) is 21.59 decibels 
(which amounts to an average error for the Legendre approximation 
0.087 times the error for the Taylor approximation).
The approximation for the sine function leads to an improvement
of 63.18 decibels (that is, the error the Legendre approximation is
0.0007 times as that for the Taylor approximation on average). Note
in Figure~\ref{fig:TayvsLagFit} that while the Taylor polynomial
approximation has maximum accuracy at the midpoint, Legendre polynomial
approximation distributes the error more uniformly throughout the
entire interval. As a result, the squared error of is smaller and
the SNR is larger in case of Legendre polynomials.

\section{Conclusions}

Our analytical and numerical results show that Legendre polynomials can 
substantially improve the speed and accuracy of function approximation 
compared to Taylor polynomials of the same order.
The fast convergence of Legendre polynomials was noted in a prior study 
\cite{li1998}, but the geometric convergence in norm has not been shown 
analytically before. The geometric convergence rate is consistent with 
the general result of Srivastava \citep{srivastava2009} on the relation 
of the generalized rate of growth of an entire function to the rate of 
uniform convergence of a polynomial approximation on an arbitrary infinite 
compact set. And, it should be noted that the geometric convergence is not 
possible in general for $L_p$ approximants if $p\neq 2$ \cite{ivanov1990}. 
Thus, approximation using Legendre polynomials can provide significant 
performance improvements in practical applications.
We also showed that Legendre polynomials has the additional advantage 
that an arbitrary function can be approximated using \emph{non-uniformly} 
spaced points on a given interval. Importantly, its accuracy of approximation 
is substantially higher than that of the Taylor expansion with the same order 
of polynomials, with a uniform error distribution across the entire interval.
Therefore, Legendre expansion, instead of Taylor expansion, should be used 
when global accuracy is important.


\end{document}